# A FLEXIBLE IMPLEMENTATION OF A MATRIX LAURENT SERIES-BASED 16-POINT FAST FOURIER AND HARTLEY TRANSFORMS

R. C. de Oliveira

Computer Engineering Department
Amazon State University
Av. Darcy Vargas, 1200, Parque 10,
69065-020, Manaus-AM, Brazil.
email: rcorrea.oliveira@gmail.com

H. M. de Oliveira, R. Campello de Souza, E.J.P. Santos

DES - Federal University of Pernambuco
Av. Ac. Hélio Ramos, s/n, CDU, 50740-530,
Recife-PE, Brazil
email: (hmo, ricardo)@ufpe.br, edval@ee.ufpe.br

## ABSTRACT

This paper describes a flexible architecture for implementing a new fast computation of the discrete Fourier and Hartley transforms, which is based on a matrix Laurent series. The device calculates the transforms based on a single bit selection operator. The hardware structure and synthesis are presented, which handled a 16-point fast transform in 65 nsec, with a Xilinx SPARTAN 3E device.

## 1. INTRODUCTION

Fourier transforms play a major role in the fields related to Signal Processing [1, 2]. The successful application of transform techniques is mainly due to the existence of the so-called fast algorithms [3]. This paper proposes a new fast algorithm and its hardware implementation, for computing the discrete Fourier (DFT) [4] and Hartley (DHT) [5] transforms of sequences of particular lengths $N$, namely those for which $N \equiv 0 \pmod{4}$. Let $N$ be the number of time-domain samples of a sequence $v = \{v_n\}$, $n = 0,1,2,...,N$-1. The DFT of $v$ is given by the sequence $V = \{V_k\}$, of length $N$, defined by Equation (1).

$$V_k := \sum_{n=0}^{N-1} v_n \exp\left(-\frac{j2\pi kn}{N}\right). \quad (1)$$

The Discrete Hartley Transform (DHT) is defined by Equation (2) [6].

$$H_k := \sum_{n=0}^{N-1} v_n \left[\cos\left(\frac{2\pi kn}{N}\right) + \sin\left(\frac{2\pi kn}{N}\right)\right]. \quad (2)$$

From equations (1) and (2), it is apparent that the DFT can be computed from the DHT by Equation (3).

$$H_k = \Re e(V_k) - \Im m(V_k). \quad (3)$$

In 1965, J.W. Cooley and J.W. Tukey introduced a revolutionary idea which later became known as the Fast Fourier Transform (FFT) [4]. The FFT is a milestone in the theory of algorithms [7-9].

In 1984, R. N. Bracewell introduced an algorithm for performing the discrete Hartley transform (DHT) [4]. With the advent of VLSI and the development of the Digital Signal Processor to implement signal processing techniques, the DFT became the most attractive tool for spectrum evaluation [10-12]. The cost reduction of DSPs and the astonishing capacity achieved by up to date processors (e.g., dozens of GFlops–Giga floating-point operations per second) [13] is turning real-time application feasible for several kind of signals. Therefore, discrete transforms became the widespread tool in spectral analysis [14]. A lucid tutorial review of fast Fourier techniques is available in [15-16]. In 2000, an algorithm based on multilayer decomposition to calculate the DFT via the DHT was introduced [17-18]. This paper proposes a flexible implementation for the FFT and the Fast Hartley Transform (FHT), using a new approach, which is derived from a matrix-based Laurent series expansion [19].

## 2. THE FFT/FHT ALGORITHM

The fast algorithm is written according to the following DFT matrix decomposition: $\Re e(M)$ and $\Im m(M)$ denotes the real and imaginary parts of the matrix $M$, equations (3) and (4), respectively [20].

$$\Re eDFT = \left\{\Re e(M_0) + \sum_{m=1}^{(N/4-1)/2} \Re e(M_m + M_{-m})\cos\frac{2\pi m}{N}\right\}$$
$$+ \left\{\sum_{k=1}^{(N/4-1)/2} \Im m(M_m - M_{-m})\sin\frac{2\pi m}{N}\right\}, \quad (4)$$

$$\Im mDFT = \left\{\Im m(M_0) + \sum_{m=1}^{(N/4-1)/2} \Im m(M_m + M_{-m})\cos\frac{2\pi m}{N}\right\}$$
$$- \left\{\sum_{k=1}^{(N/4-1)/2} \Re e(M_m - M_{-m})\sin\frac{2\pi m}{N}\right\}, \quad (5)$$

where the matrices $M_m$ are given by Equation (5).

$$M_m := \sum_{l \in C_m} (-j)^{\frac{4(l-m)}{N}} \chi_l(M). \quad (6)$$

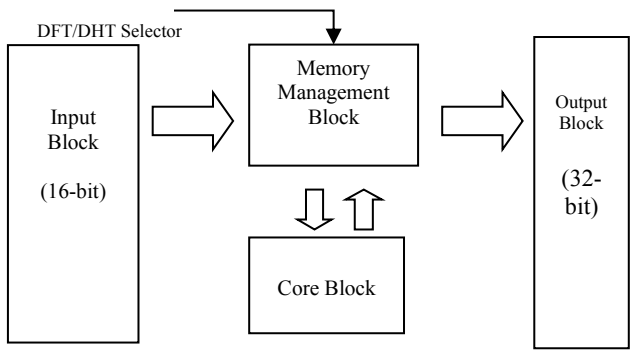

**Fig. 1 - Fast algorithm architecture FFT/FHT. The arithmetic block corresponds to the Core Block.**

The operator $\chi_l$ acts over an $N \times N$ matrix, for each $l = 0,1,2,\ldots,N-1$, yielding a new $N \times N$ binary matrix whose elements are $(\delta_{l,m})$ with $m := m_{k,n}$, where $\delta$ is the Kronecker symbol. The matrices $\Re e(M_0)$, $\Re e(M_m \pm M_{-m})$, $\Im m(M_0)$ and $\Im m(M_m \pm M_{-m})$ are then written in standard echelon form. From (4) and (5), the DHT components can be computed by Equation (3).

## 3. DESIGN METHODOLOGY AND ARCHITECTURE

The design was carried out through the steps: specification, VHDL description, behavioral simulation and synthesis.

### 3.1 Specification

The project aims at the production of a fast DFT/DHT computer, for real sequences with blocklength 16, every component of which is represented by a 16-bit word. The computations are to be made using fixed-point arithmetic (7 bits). Every component of the output sequence is represented by a 32-bit word (16 bits for its real part and 16 bits for its imaginary part). With a single bit selection, the user can choose which transform is to be computed.

### 3.2 VHDL Description

The device VHDL description was generated with the aid Matlab™ (Simulink), according to Equations (5) and (6).

### 3.3 Behavioral Simulation

With Xilins ISE tool, the VHDL description was compiled and simulated to check the output data. Syntax errors were removed at this point. The behavioral simulation was carried out by Xilinx ISE tool, and testbench files were generated.

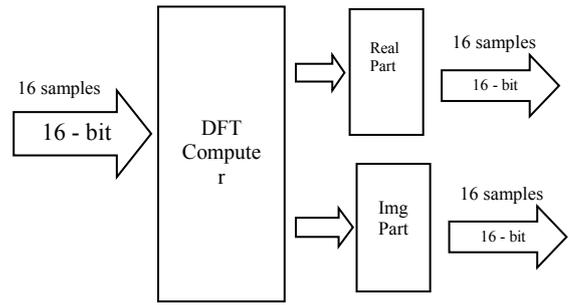

**Fig. 2 – Core Block architecture. This block is used both in DFT/DHT computation. (See Fig. 1)**

The stimulus (input values) were inserted and results colleted according to diagrams from Fig. 3 (DHT) and Fig. 4 (DFT).

### 3.4 Synthesis

At this step, the VHDL code was analysed and optimized by the synthesis tool, in order to create an efficient implementation of the device. A Register Transfer Level (RTL) scheme was generated. The construction tool generates the file to be burn-in in the chosen device, namely, the Spartan 3E, xc3s500e-5-fg320 device. In the post-synthesis simulation, the finished circuit process occurs at 65 ns, and all obtained results corroborate the previous behavioral simulation. Main characteristics of the hardware: number of slices is 1611 out of 4656 (34%), number of slice flip-flops is 656 out of 9312 (7%), number of 4 input LUTs is 2894 out of 9312 (31%), number of IOs is 56, number of bonded IOBs is 56 out of 232 (24%), number of MULT18X18SIOs is 12 out of 20 (60%), and number of GCLKs is 2 out of 24 (8%).

### 3.5. Architecture

The device architecture is based on two main blocks, the memory management and the core block, as shown in Fig. 1. The block of memory management stores the components of the input signal and the output transform.
The core block is responsible for calculating the DFT coefficients, according to Fig. 2. From these, the memory management block selects and computes the desired transform.
After storage in the memory, the device calculates the transforms based on the selected operator (DFT/DHT). Thereby, the transforms are stored in the memory, as follows. For the DFT: the first two bytes hold the real components and the last two hold the imaginary components; For the DHT: only the last 16 bits hold the real components.

## Table 1. DFT/DHT values via Behavioral Simulation Data and via MatLab<sup>TM</sup> (into brackets).

| Index ($k$) | Input Data | Output Data DFT | Output Data DHT |
|---|---|---|---|
| 0 | 0 | 56 [56.000] | 56 [56.0000] |
| 1 | 1 | 0 [0] | 0 [0] |
| 2 | 2 | -8.0000 +19.3137$j$ [-8.0000 +19.375$j$] | -27.3137 [-27.375] |
| 3 | 3 | 0 [0] | 0 [0] |
| 4 | 4 | -8.0000 + 8.0000 $j$ [-8.0000+8.0000$j$] | -16.0000 [-16.0000] |
| 5 | 5 | 0 [0] | 0 [0] |
| 6 | 6 | -8.0000 + 3.3137$j$ [-8.0000 + 3.375$j$] | -11.3137 [-11.375] |
| 7 | 7 | 0 [0] | 0 [0] |
| 8 | 0 | -8.0000 [-8.0000] | -8.0000 [-8.0000] |
| 9 | 1 | 0 [0] | 0 [0] |
| 10 | 2 | -8.0000 - 3.3137$j$ [-8.0000 - 3.375$j$] | -4.6863 [-4.625] |
| 11 | 3 | 0 [0] | 0 [0] |
| 12 | 4 | -8.0000 - 8.0000$j$ [-8.0000 - 8.0000$j$] | 0 [0] |
| 13 | 5 | 0 [0] | 0 [0] |
| 14 | 6 | -8.0000 -19.3137$j$ [-8.0000 -19.375$j$] | 11.3137 [11.375] |
| 15 | 7 | 0 [0] | 0 [0] |

The core block was derived from a Simulink<sup>TM</sup> implementation. The arithmetic complexity for this 16-point fast algorithm is 12 fixed-point multiplications (Equation 7) and 101 additions.

## 4. SIMULATION RESULTS

The implementation was simulated using Xilinx ISE tool to verify the accuracy of the device and the results were compared with the ones obtained from Simulink<sup>TM</sup>. The results are shown in Figs. 3, 4. The simulation was made with several input sequences. For the particular sequence $v$ = {0 1 2 3 4 5 6 7 0 1 2 3 4 5 6 7}, the results are shown in Table 1. The output DFT sequences of the simulation carried out by the device was exactly the same as that one obtained by Simulink<sup>TM</sup> (Column 3, Table 1). Both are slightly different than the true-DFT sequence computed by an internal MatLab<sup>TM</sup> routine (Column 3, into brackets). The quantization error ( 0.3%) is due to the use of fixed-point instead of floating-point arithmetic. Nevertheless, this error magnitude is acceptable for many applications.

## 5. CONCLUDING REMARKS

This paper presented the design and implementation on the xc3s500e-5-fg320 device, of an algorithm that is capable of fast computing the DFT/DHT coefficients of a real sequence with blocklength $N$=16. The device computes either transform based on a single bit selection. The device computed the DFT/DHT at 65 nsec, which is acceptable for applications such as audio, speech processing, biomedical signal processing and ×DSL.

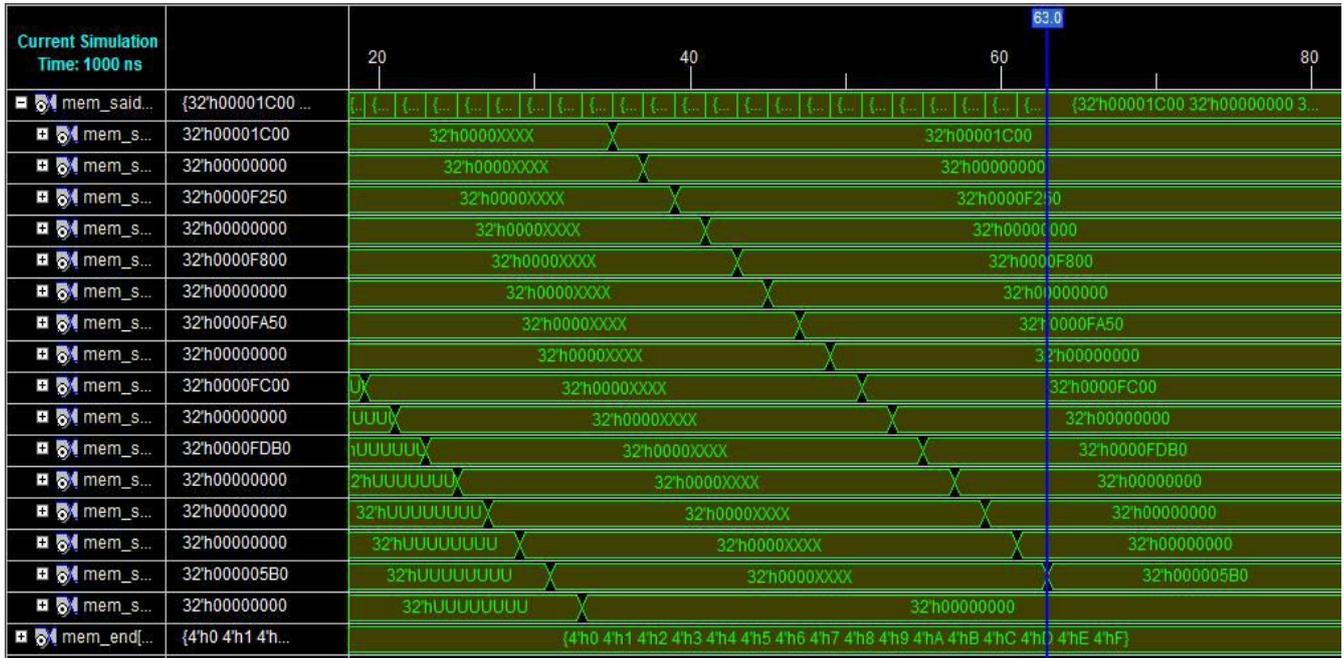

**Fig. 3** Results for behavioral simulation in operation DHT using Xilinx ISE **tool.**

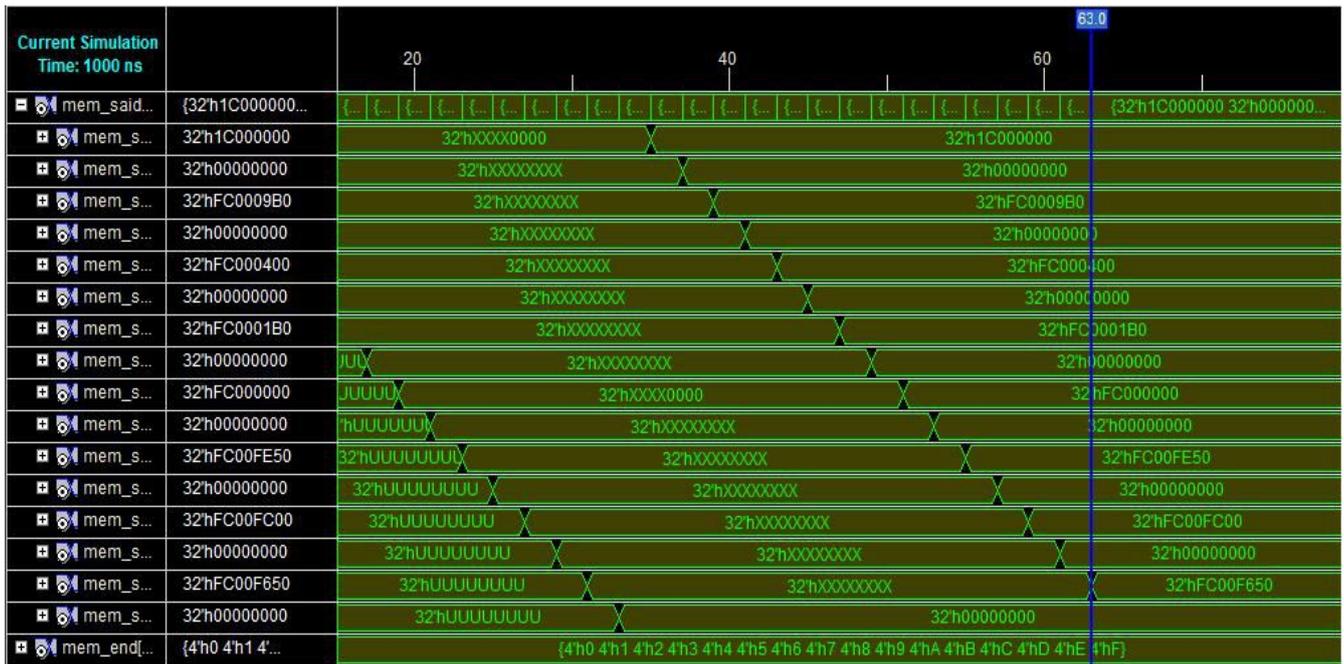

**Fig. 4** Results for behavioral simulation in operation DFT using Xilinx ISE tool.